\documentclass{article}

\usepackage{amssymb}

\newcommand{\real}{{\mathbb R}}
\newcommand{\iin}{\!\in\!}
\newcommand{\Lip}{\mbox{\rm Lip}}
\newcommand{\Lipplus}{\mbox{\rm Lip}^+}
\newcommand{\Ub}{{\bf U_b}}
\newcommand{\Ubplus}{{\bf U^+_b}}
\newcommand{\Meas}{{\bf M}}
\newcommand{\UMeas}{{\bf M_u}}
\newcommand{\UMeasplus}{{\bf M^+_u}}
\newcommand{\Mol}{{\bf Mol}}
\newcommand{\sect}[1]{\setminus_{#1}}
\newcommand{\rtrans}[2]{\rho^{#1}(#2)}
\newcommand{\orbit}[1]{\mbox{\rm orb}(#1)}
\newcommand{\clorbit}[1]{\overline{\mbox{\rm orb}(#1)}}
\newcommand{\card}[1]{| #1 |}
\newcommand{\up}{{\cal UP}}
\newcommand{\scrO}{{\cal O}}

\newcommand{\bigsep}{\mbox{\Large $|$}}
\newcommand{\conv}{\star}
\newcommand{\measprod}{\otimes}
\newcommand{\semiu}{\raisebox{1mm}{$\ast$}}
\newcommand{\wstar}{weak\raisebox{1mm}{$\ast$}\ }
\newcommand{\compln}[1]{\widehat{#1}}
\newcommand{\compactn}[1]{\overline{#1}}
\newcommand{\psm}{\Delta}
\newcommand{\qed}{\hfill $\Box$\vspace{3mm}}

\newtheorem{theorem}{Theorem}[section]
\newtheorem{lemma}[theorem]{Lemma}
\newtheorem{corollary}[theorem]{Corollary}
\newtheorem{question}{Question}

\title{Semiuniform semigroups and convolution}

\author{Jan Pachl  \\ Toronto, Ontario, Canada}
\date{November 26, 2008 (version 2)}

\begin{document}
\maketitle

\begin{abstract}
Semiuniform semigroups provide a natural setting for the
convolution of generalized finite measures on semigroups.
A semiuniform semigroup is said to be ambitable if each uniformly bounded
uniformly equicontinuous set of functions on the semigroup
is contained in an ambit.
In the convolution algebras constructed over ambitable semigroups,
topological centres have a tractable characterization.
\end{abstract}


\section{Semiuniform products and semiuniform semigroups}
    \label{section:definitions}

Functional analysis on topological and semitopological semigroups is well
developed~\cite{BJM1989}.
However, concrete semigroups often carry not only a topology but also
a natural uniform structure.
That suggests that it is worthwhile to investigate semigroups endowed
with compatible uniform structures.

An observation at the end of the author's paper~\cite{Pachl2006} points out
that the main results in that paper,
proved there for topological groups, hold more generally for semiuniform semigroups.
Moreover, the definition and basic properties of
ambitable topological groups~\cite{Pachl2008},
and their connection to topological centres in convolution algebras,
are easily generalized to semiuniform semigroups.
These generalizations are described in the current paper.

Several theorems below are obtained by simple modifications of the proofs
for topological groups
in the author's previous papers~\cite{Pachl2006}\cite{Pachl2008}.
The modified proofs are included here for the sake of completeness.

Most of the notation used here is defined in~\cite{Pachl2006}.
All uniform structures are assumed to be Hausdorff,
and all linear spaces over the field $\real$ of reals.

Following Isbell~\cite{Isbell1964}, we denote by $ X \semiu Y $
the {\em semiuniform product\/} of uniform spaces $X$ and~$Y$.
By definition,
a semigroup $X$ with a uniform structure is a {\em semiuniform semigroup\/}
if the semigroup operation $ (x,y) \mapsto xy $ is uniformly
continuous from the semiuniform product $ X \semiu X $ to $X$.
In other words, $X$ is a semiuniform semigroup if and only if
\begin{itemize}
\item the set $ \{ x \mapsto xy  \; | \; y \!\in\! X  \} $ of mappings from $X$ to $X$
is uniformly equicontinuous; and
\item for each $ x \in X $, the mapping $ y \mapsto xy  $ from $X$ to $X$
is uniformly continuous.
\end{itemize}

Hindman and Strauss~\cite{HS1998} used these two conditions in
extending a semigroup operation to the uniform compactification.

A more precise term would be a {\em right semiuniform semigroup\/},
since we may also define a {\em left semiuniform
semigroup\/} by requiring the mapping $ (x,y) \mapsto yx $ to be
uniformly continuous from $ X \semiu X $ to $X$.
However, in this paper we only study right semiuniform semigroups,
and omit the qualifier {\em right\/}.

The class of semiuniform semigroups includes several familiar classes
of semigroups with additional structure:

\begin{description}
\item[Discrete semigroups.]
Every semigroup with the discrete uniformity is a semiuniform semigroup.
\item[Topological groups.]
Every topological group with its right uniformity
is a semiuniform semigroup (\cite{Pachl2006}, Lemma~4.1).
\item[Uniform semigroups.]
Every uniform semigroup in the sense of Marxen~\cite{Marxen1973}
is a semiuniform semigroup.
\item[Semigroups of uniform endomorphisms.]
Let $Y$ be a uniform space, and let $U(Y,Y)$ be the uniform space of
uniformly continuous mappings from $Y$ to itself,
as defined by Isbell~(\cite{Isbell1964}, Ch.~III).
With the composition operation, $U(Y,Y)$ is a semiuniform semigroup.
\item[Balls in normed algebras.]
Let $A$ be a normed algebra and $r$ a real number, $ 0 < r \leq 1 $.
Let $B_r$ and $\overline{B_r}$ be the open and the closed ball with diameter $r$ in $A$.
Then $B_r$ and $\overline{B_r}$ with the algebra multiplication and the uniformity
defined by the norm in $A$ are semiuniform semigroups.
\end{description}


\section{Convolution}
    \label{section:convolution}

As in~\cite{Pachl2006},
if $X$, $Y$ and $Z$ are sets and $p$ is a mapping from $ X \times Y $ to $ Z $
then define $ \sect{x} p(x,y) $ to be the mapping
$ x \mapsto p(x,y) $ from $X$ to $Z$
and $ \sect{y} p(x,y) $ to be the mapping $ y \mapsto p(x,y) $ from $Y$ to $Z$.

When $X$ is a uniform space,
denote by $\up(X)$ the set of all uniformly continuous pseudometrics on $X$.
Recall~\cite{Pachl2006}\cite{Pachl2008} that for a uniform space $X$,
\begin{itemize}
\item
$\Ub (X)$ is the space of bounded uniformly continuous real-valued functions on $X$
with the $\sup$ norm $\| . \|$;
\item
for a pseudometric $\psm$ on $X$,
\begin{eqnarray*}
\Lip ( \psm )  & = &  \{ \; f: X \rightarrow [-1,1] \; \bigsep \;
| f(x) - f(x') | \leq \psm(x,x') \; \mbox{\rm for all} \; x,x' \iin X  \;\} \\
\Lipplus ( \psm )  & = & \{ \; f \iin \Lip(\psm) \; \bigsep \;
f(x) \geq 0 \; \mbox{\rm for all} \; x \iin X \;\}
\end{eqnarray*}
and both $\Lip(\psm)$ and $ \Lipplus ( \psm ) $
are always considered with the topology of pointwise convergence on $X$;
\item
$\Meas ( X )$ is the norm dual of $\Ub ( X )$;
\item
the {\em \wstar topology\/} on $\Meas(X)$ is the weak topology
of the duality $ \langle \, \Meas(X), \Ub(X) \, \rangle $;
\item
the {\em UEB topology\/} on $\Meas(X)$ is the topology of uniform convergence
on the sets $\Lip(\psm)$, where $\psm$ ranges
over all uniformly continuous pseudometrics on $X$;
equivalently, it is the topology of uniform convergence
on uniformly equicontinuous bounded subsets of $ \Ub(X) $.
\item
the {\em UEB uniformity\/} is the corresponding translation-invariant
uniformity on $ \Meas(X) $;
\item
when $p$ is a uniformly continuous mapping from $X$ to a uniform space $Y$,
the linear mapping $\Meas (p) : \Meas (X) \rightarrow \Meas (Y) $
is defined by $ \Meas (p) (\mu) (f) = \mu ( f \circ p ) $
for $ f \iin \Ub (Y) $;
\item
the subspace $\UMeas ( X )$ of $\Meas ( X )$ is defined as follows:
$\mu \in \UMeas ( X )$ iff $\mu$ is continuous on $\Lip(\psm)$
for each $\psm \iin \up(X)$
(or, equivalently, $\mu$ is continuous on $\Lipplus(\psm)$
for each $\psm \iin \up(X)$);
\item
if $X$ and $Y$ are uniform spaces, $\mu \iin \Meas(X)$, $\nu \iin \Meas(Y)$, then
the {\em direct product\/} $\mu \measprod \nu$ is an element of $\Meas( X \semiu Y )$
defined by $ \mu \measprod \nu (f) = \mu ( \,\sect{x} \nu ( \,\sect{y} f(x,y) )$
for $ f \iin \Ub ( X \semiu Y ) $;
\item
for $ x \in X $, $ \delta_x \in \Meas ( X ) $ is defined by
$ \delta_x (f) = f(x) $ for $ f \in \Ub ( X ) $;
\item
$\Mol(X)$ is the linear subspace of $\Meas(X)$ generated by
the set $ \{ \;\delta_x \; | \; x \in X \;\}$;
\item
the mapping $ \delta: x \mapsto \delta_x $ is a topological embedding of $X$
to $ \Meas ( X ) $ with the \wstar topology;
\item
$ \compactn{X} $,
the \wstar closure of $ \delta(X) $ in $ \Meas ( X ) $,
is a {\em uniform compactification of $X$\/};
\item
$ \compln{X} = \compactn{X} \cap \UMeas ( X )$ with the UEB uniformity
is a completion of $X$;
\item
we define
\begin{eqnarray*}
\Ubplus (X) \; & = & \; \{ \; f \in \Ub(X) \; | \; f(x) \geq 0 \; \mbox{\rm for} \;
x \in X \; \} \\
\Meas^+ (X) \; & = & \; \{ \; \mu \in \Meas(X) \; | \; \mu(f) \geq 0 \; \mbox{\rm for} \;
f \in \Ubplus (X) \; \}
\end{eqnarray*}
and similarly for $\Mol^+ (X)$ and $\UMeasplus (X)$.
\end{itemize}

Let $X$ be a semiuniform semigroup and $ \mu , \nu \iin \Meas(X) $.
Then $ \mu \measprod \nu \iin \Meas( X \semiu X)$.
Since the mapping $ m : (x,y) \mapsto xy $
is uniformly continuous from $X \semiu X$ to $X$,
the image $ \Meas(m) ( \mu \measprod \nu ) $ of $ \mu \measprod \nu $
is well defined as an element of $\Meas(X)$.
We denote
$
\mu \conv \nu = \Meas(m) ( \mu \measprod \nu )
$
and call $ \mu \conv \nu $ the {\em convolution\/} of $\mu$ and $\nu$.
Expanding the definition, we obtain
\[
\mu \conv \nu (f)  = \mu ( \sect{x} \nu ( \sect{y} f(xy ) ) )
\]
for $ \mu , \nu \in \Meas ( X ) $, $f \in \Ub ( X )$.
Thus $ \mu \conv \nu $ is the {\em evolution\/} of $\mu$ and $\nu$ in the terminology
of Pym~\cite{Pym1964}\cite{Pym1965}.
On $ \compactn{X} $, the convolution $ \conv $
coincides with the operation defined by Hindman and Strauss (\cite{HS1998}, 21.43).

More generally,
let $X$ be a semigroup {\em acting semiuniformly\/} on a uniform space $Y$.
By that we mean that $X$ is a semigroup endowed with a uniformity and
there is a uniformly continuous mapping $ m: X \semiu Y \rightarrow Y $ such that
$ m ( s , m ( s' , y ) ) = m ( s s' , y ) $
for all $ s, s' \in X $, $ y \in Y$.
In that case,
define the convolution operation from $ \Meas(X) \times \Meas(Y) $ to $ \Meas(Y) $
by $ \mu \conv \nu = \Meas(m) ( \mu \measprod \nu ) $
for $\mu \iin \Meas(X)$ and $\nu \iin \Meas(Y)$.
The definition of $ \mu \conv \nu $ for a semiuniform semigroup $X$
and $ \mu , \nu \iin \Meas(X) $
is obtained as a special case for $ X = Y $.

The proofs in this section are completely analogous to those
in section 4 in~\cite{Pachl2006},
where the same results are proved for topological groups.

\begin{theorem}
    \label{th:associative}
Let $X$ be a semiuniform semigroup acting semiuniformly on a uniform space $Y$.
Let $ \mu , \mu' \in \Meas(X) $, $ \nu \in \Meas(Y) $.
Then
$
( \mu \conv \mu' ) \conv \nu = \mu \conv ( \mu' \conv \nu )
$.
\end{theorem}

\noindent
{\bf Proof.} Apply Lemma 3.7 in~\cite{Pachl2006}.
\qed

\begin{theorem}
    \label{th:identities}
Let $X$ be a semigroup acting semiuniformly on a uniform space $Y$.
Let $ r \in \real $,
$ \mu , \mu' \in \Meas(X) $, $ \nu, \nu' \in \Meas(Y) $.
Then
\begin{eqnarray*}
& & ( r \mu ) \conv \nu = \mu \conv ( r \nu ) = r ( \mu \conv \nu )     \\
& & ( \mu + \mu' ) \conv \nu = ( \mu \conv \nu ) + ( \mu' \conv \nu )   \\
& & \mu \conv ( \nu + \nu' ) = ( \mu \conv \nu ) + ( \mu \conv \nu' )
\end{eqnarray*}
\end{theorem}

\noindent
{\bf Proof.}
Apply Lemma 3.6 in~\cite{Pachl2006}.
\qed

\begin{theorem}
    \label{th:closure}
Let $X$ be a semigroup acting semiuniformly on a uniform space $Y$.
Let $\Phi$ be one of the functors $\Meas^+$, $\Mol$, $\Mol^+$, $\UMeas$, $\UMeasplus$,
$\sect{X} \compactn{X}$ and $\sect{X} \compln{X}$.
If $\mu \iin \Phi(X)$ and $\nu \iin \Phi(Y)$ then $ \mu \conv \nu \iin \Phi(Y) $.
\end{theorem}

\noindent
{\bf Proof.}
Apply Lemma~3.5, Theorem~3.13 and Lemma~6.1 in~\cite{Pachl2006}.
\qed

\begin{theorem}
    \label{th:convolution-topology}
Let $X$ be a semigroup acting semiuniformly on a uniform space $Y$.
\begin{enumerate}
\item
    \label{th:convolution-topology:part1}
$ \| \mu \conv \nu \| \leq \| \mu \| . \| \nu \| $
for any $ \mu \in \Meas(X) $, $ \nu \in \Meas(Y) $.
\item
    \label{th:convolution-topology:part2}
Let $ B \subseteq \Meas(Y) $ be a set bounded in the $\| . \|$ norm,
$ \mu_0 \in \UMeas(X) $, and $ \nu_0 \in B $.
When $ \Meas(X)$ and $B$ are endowed with their UEB topology,
the mapping $ ( \mu, \nu ) \mapsto \mu \conv \nu $
from $ \Meas(X) \times B $ to $ \Meas(Y) $
is jointly continuous at $ ( \mu_0 , \nu_0 ) $.
\item
    \label{th:convolution-topology:part3}
The mapping $ ( \mu, \nu ) \mapsto \mu \conv \nu $
is jointly sequentially continuous
from $ \UMeas(X) \times \UMeas(Y) $ to $ \UMeas(Y) $
when $ \UMeas(X)$ and $\UMeas(Y)$ are endowed with their \wstar topology.
\item
    \label{th:convolution-topology:part4}
Let $ \mu_0 \in \UMeasplus(X) $ and $ \nu_0 \in \Meas^+ (Y) $.
When $ \Meas^+(X) $ and  $ \Meas^+(Y) $ are endowed with their \wstar topology,
the mapping $ ( \mu, \nu ) \mapsto \mu \conv \nu $
from $ \Meas^+(X) \times \Meas^+(Y) $ to $ \Meas^+ ( Y ) $
is jointly continuous at $ ( \mu_0 , \nu_0 ) $.
\item
    \label{th:convolution-topology:part5}
If $ \mu \in \UMeas(X) $ then the mapping $ \nu \mapsto \mu \conv \nu $
from $ \Meas(Y) $ to itself is \wstar continuous.
\end{enumerate}
\end{theorem}

\noindent
{\bf Proof.}
All statements follow from results in section~3 of~\cite{Pachl2006}:
Part~\ref{th:convolution-topology:part1}
from Lemma~3.4,
Part~\ref{th:convolution-topology:part2}
from Theorem~3.8,
part~\ref{th:convolution-topology:part3}
from Corollary~3.9,
part~\ref{th:convolution-topology:part4}
from Theorem~3.10
and part~\ref{th:convolution-topology:part5}
from Theorem~3.12.
\qed

As a corollary we obtain that for any semiuniform semigroup $X$
the spaces $ \Meas(X) $ and $ \UMeas(X) $ with
the operations $ \conv $ and $ + $ and the norm $ \| . \| $ are Banach algebras.

For the special case of topological groups,
it is noted in section~4 of~\cite{Pachl2006} that algebraic identities are
inherited from $\Mol(X)$ to $\UMeas(X)$.
In view of Theorem~\ref{th:convolution-topology},
the same is true for semiuniform semigroups.
In particular, if a semiuniform semigroup $X$ is commutative then so is $\UMeas(X)$.


\section{Ambitable semigroups and topological centres}
    \label{section:centre}

The following definitions and questions are a straightforward generalization of those for
ambitable topological groups~\cite{Pachl2008}.

Let $X$ be a semigroup, $f$ a real-valued function on $X$ and $ x \in X $.
The {\em right translation\/} of $f$ by $x$ is the function
$\rtrans{x}{f} = \sect{z} f(zx) $.
The set
$ \orbit{f} = \{ \rtrans{x}{f} \; | \; x \in X \} $
is the {\em right orbit\/} of $f$, and
$ \clorbit{f} $ is
the closure of $ \orbit{f} $ in the product space $ \real^X $.

Say that a semiuniform semigroup $X$ is {\em ambitable\/} if for every
$ \psm \iin \up(X) $
there exists $ f \in \Ub(X) $ such that
$ \Lipplus(\psm) \subseteq \clorbit{f} $.

\begin{question}
\label{question:ambitable}
Which semiuniform semigroups are ambitable?
\end{question}

This question is motivated by
Questions~\ref{question:centreM} and~\ref{question:centreC} below,
in the same way as for topological groups~\cite{Pachl2008}.

If $S$ is a semigroup with a topology, define its {\em (right) topological centre\/} by
\[
\Lambda (S) =  \{ \; x \iin S \;
| \;
\mbox{\rm the mapping} \;\; y \mapsto xy \;\;
\mbox{\rm is continuous on} \; S \; \} \;.
\]

Now let $X$ be a semiuniform semigroup.
Consider the semigroup $\Meas(X)$ with the convolution operation $\conv$ and the
\wstar topology, and its subsemigroup $\compactn{X}$ (uniform compactification of $X$).
By part~\ref{th:convolution-topology:part5} of~Theorem~\ref{th:convolution-topology},
$ \UMeas(X) \subseteq \Lambda( \Meas(X) ) $,
and therefore also $ \compln{ X } \subseteq \Lambda ( \compactn{X} ) $.

The results of Lau~\cite{Lau1986}, Lau and Pym~\cite{Lau-Pym1995},
Lashkarizadeh Bami~\cite{LBami2000},
Ferri and Neufang~\cite{Ferri-Neufang2007},
and Dales, Lau and Strauss~\cite{DLS2008},
among others, lead to the following questions.

\begin{question}
    \label{question:centreM}
Which semiuniform semigroups $X$ satisfy $ \UMeas(X) = \Lambda( \Meas(X) ) $ ?
\end{question}

\begin{question}
    \label{question:centreC}
Which semiuniform semigroups $X$ satisfy $ \compln{ X } = \Lambda ( \compactn{X} ) $ ?
\end{question}

If $X$ is precompact then $ \UMeas(X) = \Meas(X) $ and $ \compln{X} = \compactn{X} $,
and therefore $ \UMeas(X) = \Lambda( \Meas(X) ) $ and
$ \compln{ X } = \Lambda ( \compactn{X} ) $.

Next we repeat the proofs in section~5 of~\cite{Pachl2008} to derive positive answers to
Questions~\ref{question:centreM} and~\ref{question:centreC} for ambitable semigroups.

\begin{lemma}
    \label{lemma:halfconv}
{\rm (cf. 5.1 in~\cite{Pachl2008})}
Let $X$ be a semiuniform semigroup and $ f \in \Ub(X) $.
\begin{enumerate}
\item
The mapping $ \varphi : \nu \mapsto \sect{x} \nu ( \sect{y} f (xy) ) $ is continuous
from $ \compactn{X} $  to the product space $ \real^X $.
\item
$ \varphi ( \compactn{X} ) = \clorbit{f} $.
\end{enumerate}
\end{lemma}

\noindent
{\bf Proof.}
1. As noted above, $ \delta_x \in \Lambda ( \compactn{X} ) $ for each $ x \in X $,
and thus the mapping
$ \nu \mapsto \delta_x \conv \nu $ is \wstar continuous from
$ \compactn{X} $ to itself.
Since $ \delta_x \conv \nu (f)  = \nu ( \sect{y} f(xy) )$,
this means that the mapping
$ \nu \mapsto \nu ( \sect{y} f(xy) ) $ from
$ \compactn{X} $ to $\real$ is continuous for each $ x \in X $,
and therefore
the mapping $ \nu \mapsto \sect{x} \nu ( \sect{y} f (xy) ) $ is continuous
from $ \compactn{X} $ to $ \real^X $.

2. $ \varphi ( \delta_x ) = \rtrans{x}{f} $ for all $ x \in X $,
and therefore $ \varphi ( \delta ( X ) ) = \orbit{f} $.
The mapping $ \varphi $ is continuous by part 1, $ \compactn{X} $ is compact, and
$ \delta(X) $ is dense in $ \compactn{X} $.
It follows that $ \varphi ( \compactn{X} ) = \clorbit{f} $.
\qed

\begin{lemma}
    \label{lemma:orbitcontinuity}
{\rm (cf. 5.2 in~\cite{Pachl2008})}
Let $X$ be any semiuniform semigroup, $ \mu \in \Meas(X) $ and $ f \in \Ub(X) $.
If the mapping $ \nu \mapsto \mu \conv \nu $ from $ \compactn{X} $ to $ \Meas(X) $
is \wstar continuous then $ \mu $ is continuous on $ \clorbit{f} $.
\end{lemma}

\noindent
{\bf Proof.}
As in Lemma~\ref{lemma:halfconv}, define
$ \varphi ( \nu ) =  \sect{x} \nu ( \sect{y} f (xy) ) $
for $ \nu \in \compactn{X} $.

\begin{picture}(300,90)
\thicklines
\put(54,60){\makebox(30,20){$\compactn{X}$}}
\put(140,60){\makebox(40,20){$\clorbit{f}$}}
\put(140,0){\makebox(40,20){$\real$}}

\put(80,70){\vector(1,0){60}}
\put(160,60){\vector(0,-1){40}}
\put(80,60){\vector(3,-2){66}}

\put(95,66){\makebox(30,20){$\varphi$}}
\put(160,30){\makebox(20,20){$\mu$}}
\put(50,20){\makebox(80,20){$\nu \mapsto \mu \conv \nu (f) $}}
\end{picture}

By the definition of convolution,
$ \mu \conv \nu (f)  = \mu ( \sect{x} \nu ( \sect{y} f(xy) ) )
= \mu ( \varphi ( \nu ) ) $.
Thus $ \mu \circ \varphi $
is continuous from $ \compactn{X} $ to $\real$.

By Lemma~\ref{lemma:halfconv},
$ \varphi $ is continuous
from $ \compactn{X} $ to $ \clorbit{f} $,
and $ \varphi ( \compactn{X} ) = \clorbit{f} $.
Since $ \compactn{X} $ is compact, it follows that
$\mu$ is continuous on $\clorbit{f}$.
\qed

\begin{theorem}
    \label{th:ambitableunif}
{\rm (cf. 5.3 in~\cite{Pachl2008})}
If $X$ is an ambitable semiuniform semigroup,
$ S \subseteq \Meas(X) $,
and $S$ with the $\conv$ operation is a semigroup
such that $ \compactn{X} \subseteq S $,
then $ \Lambda (S) = \UMeas(X) \cap S $.
\end{theorem}

\noindent
{\bf Proof.}
As is noted above, $ \UMeas(X) \subseteq \Lambda ( \Meas(X) )$.
Therefore $ \UMeas(X) \cap S \subseteq \Lambda (S) $ for every semigroup
$ S \subseteq \Meas(X) $.

To prove the opposite inclusion,
take any $ \mu \iin \Lambda (S) $ and any $ \psm \iin \up(X) $.
Since $ \compactn{X} \subseteq S $,
the mapping $ \nu \mapsto \mu \conv \nu $ from $ \compactn{X} $ to $ \Meas(X) $
is \wstar continuous by the definition of $ \Lambda (S) $.
Since $ X $ is ambitable, $ \Lipplus(\psm) \subseteq \clorbit{f} $
for some $ f \in \Ub(X) $.
By Lemma~\ref{lemma:orbitcontinuity},
$\mu$ is continuous on  $ \clorbit{f} $ and therefore
also on $ \Lipplus(\psm) $.
Thus $ \mu \in \UMeas(X) $.
\qed

\begin{corollary}
    \label{corollary:twosemigroups}
{\rm (cf. 5.4 in~\cite{Pachl2008})}
If $X$ is an ambitable semiuniform semigroup then
$ \UMeas(X) = \Lambda ( \Meas(X) ) $
and
$ \compln{ X } = \Lambda ( \compactn{X} ) $.
\end{corollary}

\noindent
{\bf Proof.}
Apply \ref{th:ambitableunif} with $ S = \Meas(X) $ and with $ S = \compactn{X}$.
\qed


\section{Discrete semigroups}
    \label{section:discrete}

If $X$ is a discrete uniform space then
$ \Ub(X) = \ell^\infty (X)$,
$ \compln{X} = X $,
$\compactn{X} = \beta X $ is the Stone-\v{C}ech compactification of $X$,
and $\UMeas(X) = \ell^1 (X)$.
For discrete semigroups,
Questions~\ref{question:centreM} and~\ref{question:centreC}
of the previous section
have been extensively studied,
often in the following equivalent form:
\begin{itemize}
\item
For which semigroups $X$ is $ \ell^1 (X) $ left strongly Arens irregular?
(\cite{DLS2008},~2.24)
\item
Which semigroups are left strongly Arens irregular?
(\cite{DLS2008},~6.11)
\end{itemize}
Dales, Lau and Strauss~\cite{DLS2008}
give an example of a countable infinite abelian semigroup $X$ for which
$ \ell^1 (X) $ is not left strongly Arens irregular but $X$ is (loc.~cit., 12.21).
In other words,
$ \UMeas(X) \neq \Lambda( \Meas(X) ) $ and
$ \compln{ X } = \Lambda ( \compactn{X} ) $.
By Corollary~\ref{corollary:twosemigroups}, such $X$ is not ambitable.

This section offers a partial answer to Question~\ref{question:ambitable} for
discrete semigroups.
The approach is similar to that in section~4 of~\cite{Pachl2008},
which in turn is a variant of the factorization method in~\cite{Ferri-Neufang2007}.

When $P$ and $R$ are subsets of a semigroup $X$, write
\[
R^{-1} P  =  \{ \; x \in X \; | \; r x \in P \;\;\mbox{\rm for some}\;\; r \in R \; \}
\]

For any semigroup $X$, define the metric $\psm_X$ on $X$ by
$ \psm_X ( x , y ) = 1 $ for $ x \neq y $.
Then $ \Lipplus ( \psm_X ) = [0,1]^X $, and
$ \Lipplus ( \psm_X ) \supseteq \Lipplus ( \psm ) $
for any pseudometric $\psm$ on $X$.

The cardinality of a set $A$ is denoted $\card{A}$.
When $X$ is an infinite semigroup, consider the following two properties:
\begin{itemize}
\item[(1)]
If $F \subseteq X$ is finite then
$
\card{ \;\{ z \in X \; | \;xz \neq yz  \;\;\mbox{\rm for all}\;\;
         x,y \iin F, x \neq y \;\} \; } = \card{X}
$.
\item[(2)]
If $ P \subseteq X $, $\card{P} < \card{X}$ and $ x \iin X $ then
$ \card{{\{ x \} }^{-1} P} < \card{X} $.
\end{itemize}

Note that (1) holds in any infinite right cancellative semigroup
(in fact in any infinite near right cancellative semigroup as defined in~\cite{DLS2008}),
and (2) holds in any infinite weakly left cancellative semigroup.

The following equivalent form of (2) will be used below:
\begin{itemize}
\item[(2a)]
If $ F, P \subseteq X $, $F$ is finite and $\card{P} < \card{X}$ then
$ \card{F^{-1} P} < \card{X} $.
\end{itemize}

\begin{lemma}
    \label{lemma:finsets}
Let $X$ be an infinite semigroup satisfying (1) and (2).
Let $A$ be a set such that $\card{A} \leq \card{X}$, and
for each $ \alpha \iin A $ let $ F_\alpha $ be a non-empty finite subset of $X$.
Then there exist elements $ x_\alpha \iin X $ for $ \alpha \iin A $ such that
\begin{itemize}
\item[(i)]
the mapping $ x \mapsto x x_\alpha $
is injective on $ F_\alpha $ for every $ \alpha \iin A $ ;
\item[(ii)]
$ F_\alpha x_\alpha \cap F_\beta x_\beta = \emptyset $
for all $ \alpha , \beta \iin A $, $ \alpha \neq \beta $.
\end{itemize}
\end{lemma}

\noindent
{\bf Proof.}
Without loss of generality, assume that
$A$ is the set of ordinals of cardinality $< \card{X}$.
The construction of $ x_\alpha $ proceeds by transfinite induction.
For $ \gamma \in A $, let $ S( \gamma ) $ be the statement
``for all $ \alpha \leq \gamma $ there exist elements $ x_\alpha \in X $ such that
(i) and (ii) hold for $ \alpha, \beta \leq \gamma \;$.''

To prove $S(0)$, observe that (1) implies
\[
\card{ \;\{ z \in X \; | \;xz \neq yz  \;\;\mbox{\rm for all}\;\;
         x,y \iin F_0 , \; x \neq y \;\} \; } = \card{X}
\]
and therefore there is $ x_0 \iin X $ such that the mapping
$ x \mapsto x x_0 $ is injective on $F_0$.

Now assume that $\gamma \in A$, $\gamma > 0$,
and $S(\gamma')$ is true for all $ \gamma' < \gamma $.
We want to prove $S(\gamma)$.

As in the base step, we get
\[
\card{ \;\{ z \in X \; | \;xz \neq yz  \;\;\mbox{\rm for all}\;\;
         x,y \iin F_\gamma , \; x \neq y \;\} \; } = \card{X} \; .
\]
Let $ P = \bigcup_{\alpha < \gamma} F_\alpha x_\alpha $.
Since $ \card{P} < \card{X} $, from (2a) we get
$ \card{F_\gamma^{-1} P} < \card{X} $.
Thus there exists $ x_\gamma \iin X \setminus F_\gamma^{-1} P $ such that
the mapping $ x \mapsto x x_\gamma $ is injective on $F_\gamma$.
Since $ x_\gamma \not\in F_\gamma^{-1} P $, we have
$ F_\alpha x_\alpha \cap F_\gamma x_\gamma = \emptyset $
for all $ \alpha < \gamma $.
\qed

\begin{theorem}
    \label{th:discrete}
Every infinite discrete semigroup $X$ satisfying (1) and (2) is ambitable.
\end{theorem}

Thus, in particular, every infinite discrete group is ambitable.
That follows also from the general results in~\cite{Pachl2008}.

\noindent
{\bf Proof.}
The topology of the compact space $ \Lipplus ( \psm_X ) = [0,1]^X $
has an open basis $\scrO$ such that $\card{\scrO} = \card{X}$
and every set $ U \iin \scrO $ is a basic neighbourhood of the form
\[
U = \{ \; f \in \Lipplus (\psm_X ) \; \bigsep \;\; | f(x) - h_U (x) | < \varepsilon_U
 \; \mbox{\rm for} \; x \in F_U \; \}
\]
where $ F_U \subseteq X $ is a finite set, $h_U \iin \Lipplus(\psm_X ) $, and
$ \varepsilon_U > 0 $.

By Lemma~\ref{lemma:finsets} with $\scrO$ in place of $A$,
there are $ x_U \iin X $ for $ U \iin \scrO $ such that
\begin{itemize}
\item[(i)]
the mapping $ x \mapsto x x_U $
is injective on $ F_U $ for every $ U \iin \scrO $ ;
\item[(ii)]
$ F_U x_U \cap F_V x_V = \emptyset $
for all $ U , V \iin \scrO $, $ U \neq V $.
\end{itemize}

Define the function $ f : X \rightarrow [0,1] $ by
\[
f( y ) = \left\{
           \begin{array}{ll}
             h_U (x) & \mbox{if $y=x x_U$, $x \iin F_U$} \\[4pt]
             0       & \mbox{if $\displaystyle y \not\in \bigcup_{U \in \,\scrO} F_U x_U $}
           \end{array}
         \right.
\]
In view (i) and (ii), $f$ is well defined.
Since $ \rtrans{x_U}{f} (x) = h_U (x) $ for $x \iin F_U$,
$ \orbit{f}$ intersects every set in $\scrO$.
Thus $ \orbit{f}$ is dense in $ \Lipplus ( \psm_X ) = [0,1]^X $.
\qed

By Corollary~\ref{corollary:twosemigroups}
we obtain the following partial answers to Questions~\ref{question:centreM}
and~\ref{question:centreC}.

\begin{corollary}
If $X$ is an infinite semigroup satisfying (1) and (2) then
$ \ell^1 (X) = \Lambda( \ell^\infty (X) ) $
and $ X = \Lambda( \compactn{X} ) $.
\end{corollary}

In other words, if $X$ is an infinite semigroup satisfying (1) and (2)
then $ \ell^1 (X) $ and $X$ are left strongly Arens irregular.
This is a variant of Corollary~12.16 in~\cite{DLS2008}.


\end{document}